\documentclass[11pt,a4paper]{article}

\usepackage[headinclude]{typearea}
\usepackage[english]{babel}           
\usepackage[T1]{fontenc}
\usepackage{scrtime}
\usepackage{amsmath,amsthm,amsfonts,amssymb}
\usepackage {color}                                  
\usepackage {pifont}                                 
\usepackage {multicol,multirow}                               
\usepackage[numbers]{natbib}                                                     
\usepackage[normalem]{ulem}                                                     
\usepackage [scaled=.90]{helvet}                     
\usepackage {courier}                                
\usepackage {graphicx}                               
\usepackage {array}                                  
\usepackage {longtable}                              
\usepackage{fancyvrb}
\usepackage{enumerate} 
\usepackage{ifpdf}
\usepackage{caption}
\usepackage{mathrsfs}

\usepackage{setspace}

\usepackage{marvosym}

\usepackage{mathtools}
\mathtoolsset{showonlyrefs} 

\newtheorem{theorem}{Theorem}[section]

\theoremstyle{definition}

\newtheorem{definition}[theorem]{Definition}
\newtheorem{example}[theorem]{Example}

\newtheorem{algorithm}[theorem]{Algorithm}

\newcommand{\exclude}[1]{}

\newcommand{\E}{{\mathbb{E}}}

\newcommand{\N}{{\mathbb{N}}}
\renewcommand{\P}{{\mathbb{P}}}

\newcommand{\R}{{\mathbb{R}}}

\definecolor{darkgreen}{rgb}{0,0.5,0}
\definecolor{lightgreen}{rgb}{0.5,0.9,0.5}
\definecolor{magenta}{rgb}{0.75,0,0.25}

\definecolor{violet}{rgb}{0.25,0,0.75}

\renewcommand{\P}{{\mathbb P}}

\newcommand{\cF}{{\cal F}}

\newcommand{\cN}{{\cal N}}

\newcommand{\be}{\begin{equation}}
\newcommand{\ee}{\end{equation}}
\newcommand{\bea}{\begin{eqnarray}}
\newcommand{\eea}{\end{eqnarray}}
\newcommand{\beast}{\begin{eqnarray*}}
\newcommand{\eeast}{\end{eqnarray*}}
\newcommand{\bproof}{\begin{proof}}
\newcommand{\eproof}{\end{proof}}

\newcommand{\sign}{\operatorname{sign}}

\newcommand{\appxh}{\hat X}

\newcommand{\es}{{\underline{s}}}

\title{Stochastic differential equations with irregular coefficients:~mind the gap!}

\author{Michaela Sz\"olgyenyi\footnote{Department of Statistics, University of Klagenfurt, Universit\"atsstra\ss{}e 65--67, 9020 Klagenfurt, Austria.
michaela.szoelgyenyi@aau.at} \footnote{M.~Sz\"olgyenyi is supported by the Austrian Science Fund (FWF):~DOC 78.}}
\date{January 2021}

\begin{document}

\maketitle


Numerical methods for stochastic differential equations with non-globally Lipschitz coefficients are currently studied intensively. This article gives an overview of our work for the case that the drift coefficient is potentially discontinuous complemented by other important results in this area. To make the topic accessible to a broad audience, we begin with a heuristic on SDEs and a motivation.


\section{SDEs in a nutshell}
\label{nutshell}

A dynamical system is usually described as a solution to a deterministic (ordinary) differential equation (ODE)
\begin{align*}
x\colon [0,T]\to \R,\qquad
\frac{d}{dt} x(t)=\mu(x(t)),\qquad
x(0)=1.
\end{align*}
If we allow for random disturbance $W$ scaled by a coefficient $\sigma$ in the dynamics, this becomes
\begin{align*}
\frac{d}{dt} X(t)=\mu(X(t)) + \sigma(X(t))\frac{d}{dt} W_t,\qquad
X(0)=1.
\end{align*}
An issue here is that we cannot formally differentiate the (Wiener) noise process $W$. Hence, in stochastic analysis the \textit{stochastic differential equation (SDE)} is only a formal notation for the corresponding integral equation, taking the form
\begin{align*}
dX_t=\mu(X_t) dt + \sigma(X_t)dW_t,\qquad
X(0)=1.
\end{align*}
A solution to a stochastic differential equation is a stochastic process\footnote{A \textit{stochastic process} is a function $X\colon [0,T]\times\Omega\to\R^d$ where for all $t\in[0,T]$, $X(\cdot,t)$ is measurable. Thereby $(\Omega,\cF,\P)$ is a probability space.}.
We illustrate the difference between the deterministic and the stochastic case for $\mu\equiv-1.5$ and constant $\sigma$.
\begin{center}
\includegraphics[width=0.45\textwidth]{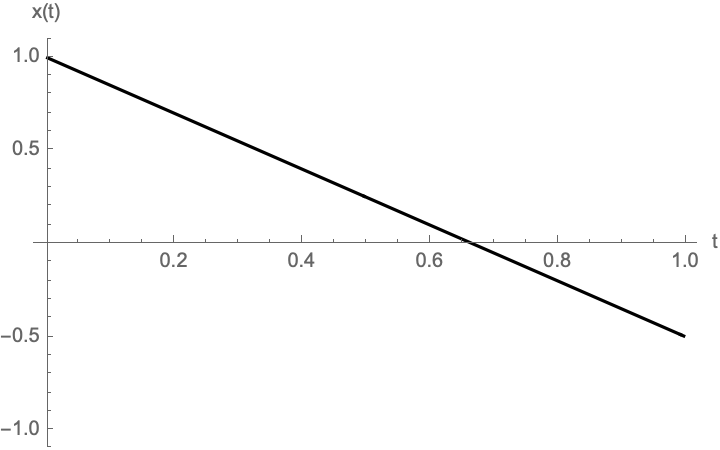}
\includegraphics[width=0.45\textwidth]{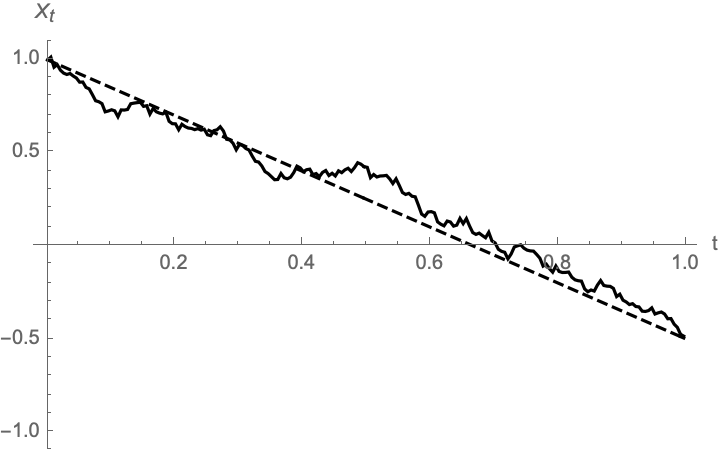}
\end{center}
The graph on the left-hand side shows the solution to the ODE, the graph on the right-hand side displays the solution to the ODE and a solution path\footnote{For every $\omega\in\Omega$, $X(\omega,\cdot)$ is called a \textit{path} of the stochastic process $X$.} of the SDE. We observe that the solution path ``zick-zacks'' randomly around the trend.

\section{Numerical methods for SDEs}
\label{numSDE}

Our objects of study are SDEs of the form
\begin{align}\label{SDE}
dX_t=\mu(X_t) dt + \sigma(X_t) dW_t, \qquad X_0=\xi,
\end{align}
where $\xi\in \R^d$, $\mu\colon\R^d\to \R^d$, $\sigma\colon\R^d\to \R^{d\times d}$, and
$W$ is a $d$-dimensional Brownian motion.
The coefficient $\mu$ is called \textit{drift}, $\sigma$ is called \textit{diffusion}.

Our example in Section \ref{nutshell} was solved explicitly. However, as in the case of ODEs, for solving most of the SDEs we have to resort to numerics. The simplest numerical method is the \textit{Euler--Maruyama method}. This approximation scheme is very similar to the Euler method for ODEs, only the noise term is added. We denote the Euler--Maruyama approximation of the solution to SDE \eqref{SDE} by $X^{\delta}$; $\delta$ is the maximal step size.
\begin{algorithm} The Euler--Maruyama method is defined as follows.
\begin{itemize}
\item Choose a time grid $0=t_0<t_1<\dots<t_n=T$ and set $\delta:=\max\{t_{k+1}-t_k\}$.
\item Start at time 0:~$X^{\delta}_{t_0}=\xi$
\item Choose $\Delta_{k+1} \sim \cN(0,t_{k+1}-t_k)$ and set
\begin{align*}
 X^{\delta}_{t_{k+1}}= X^{\delta}_{t_k}+\mu( X^{\delta}_{t_k}) \cdot (t_{k+1}-t_k) + \sigma( X^{\delta}_{t_k}) \cdot \Delta_{k+1}.
\end{align*}
\end{itemize}
\end{algorithm}

As simple and explicit schemes are desirable, we stick to the Euler--Maruyama method here. The question is now:~does this method work? The following theorem by Maruyama gives a first answer.
\begin{theorem}
 If $\mu$ and $\sigma$ are globally Lipschitz continuous, then the Euler--Maruyama method\footnote{In this paper convergence results are always given for the time-continuous versions of the schemes.} has $L^2$-convergence order $1/2$, that is $$\left(\E\!\left[ \|X_t- X^{\delta}_t \|^2 \right]\right)^{1/2}
 \le c \cdot \delta^{1/2}.$$
\end{theorem}

\subsection*{SDEs with irregular coefficients}

Now we may ask ourselves:~are we happy with this? If the coefficients that appear in your application satisfy global Lipschitz conditions, then your answer might be yes!
However in many applications the coefficients are not globally Lipschitz. In that case we speak of \textit{irregular coefficients}. Here we are particularly concerned with SDEs with discontinuous drift. This situation frequently appears, for example, in optimal control. Think of a light switch:~the light can either be on or off, influencing the drift of some underlying stochastic system. Similar discontinuities appear in energy market models, or in insurance mathematics and mathematical finance, where switching on the light is for example replaced by paying dividends to shareholders. This motivates studying SDEs with irregular coefficients from an application point of view.
But -- to be honest -- the most important motivation for many of us is mathematical curiosity.
When studying SDEs with discontinuous drift, the first question that might come to your mind is

\paragraph*{Do discontinuities matter?}

First, let us think of the ODE we solved in Section \ref{nutshell}, where we chose $\mu\equiv-1.5$ (plotted on the left-hand side below). Then, introduce a discontinuity in the point 0 (plotted on the right-hand side below).
\begin{center}
\includegraphics[width=0.45\textwidth]{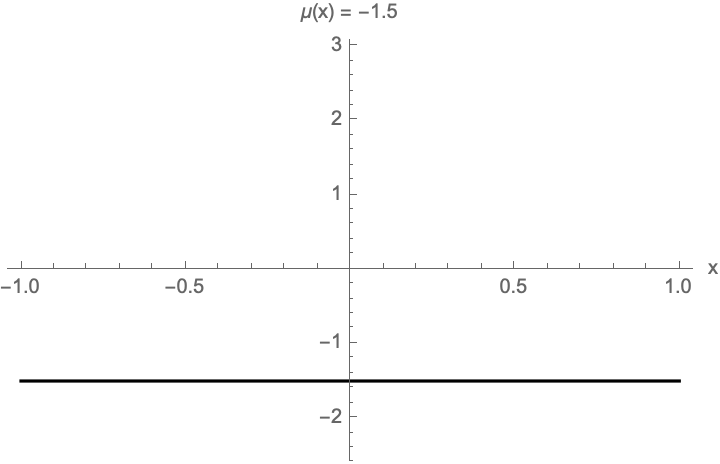}
\includegraphics[width=0.45\textwidth]{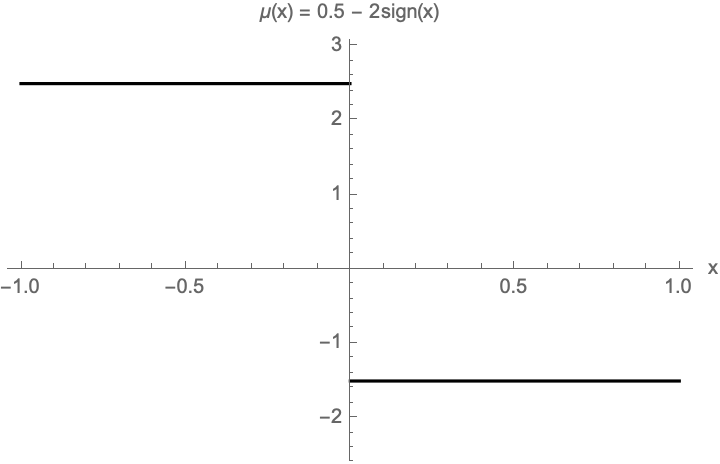}
\end{center}
The ODE corresponding to the drift on the left-hand side is explicitly solvable as we saw in Section \ref{nutshell}.
The ODE corresponding to the drift on the right-hand side,
\begin{align}\label{ode-nosol}
 \frac{d}{dt}x(t)= 0.5-2 \sign(x(t)),\qquad
 x_0=1,
\end{align}
does not admit a solution. Heuristically speaking, the reason for this is that whenever the trajectory hits zero, it is simultaneously pushed to the positive and to the negative and cannot proceed. This is illustrated in following figure\footnote{The figure was kindly provided by Gunther Leobacher.}.
\begin{center}
\includegraphics[scale=1]{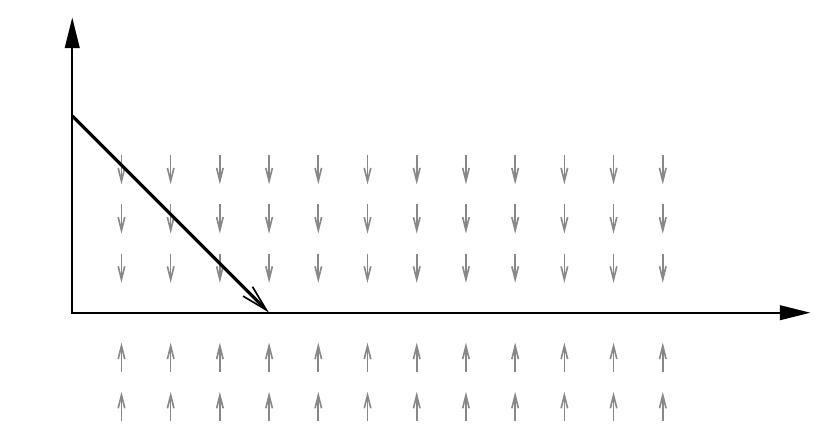}
\end{center}
This answers our question to the extent, that discontinuities do indeed matter for the question of existence of solutions of ODEs.
It also shows that one has to be careful with applying numerical solvers. While the Euler method would actually lead to a result for ODE \eqref{ode-nosol}, this result would be nonsensical.

Now let us add noise again:
\begin{align*}
 dX_t =  (0.5-2 \sign(X_t) )dt + \sigma dW_t,\qquad
  X_0=1.
\end{align*}
In this case a solution exists. In fact, the Brownian noise has a regularising effect. The noise constantly pushes the trajectory away from the point of discontinuity 0, see below. \begin{center}
\includegraphics[width=0.45\textwidth]{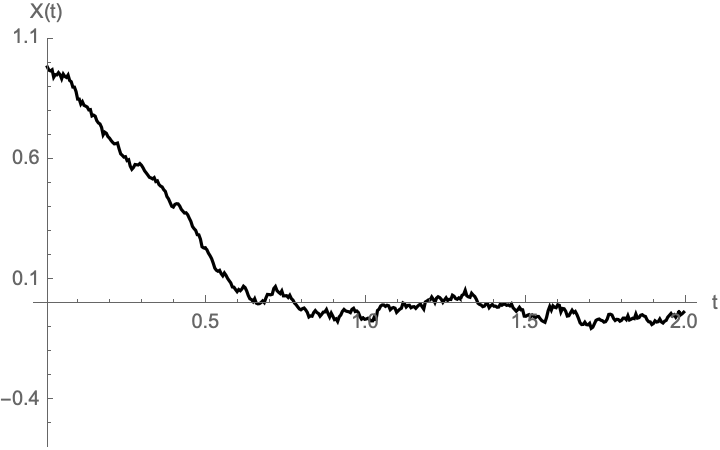}
\end{center}
Existence and uniqueness of solutions to SDEs with discontinuous drift is true under very general conditions and has been studied in various setups, see, e.g., \cite{Zvonkin1974,Veretennikov1981,Veretennikov1982,Veretennikov1984,krylov2005,zhang2013,LST15,SS16,LS16,xie2016,LS17,xie2020,PS20,leobacher2020,PSX20,soenmez2020}.\\

Now that existence and uniqueness of solutions is settled, let us turn to numerics. Can irregularities in the coefficients matter when it comes to numerics?
This leads us to a number of papers published between 2015 and 2018 where the authors construct explicit examples of SDEs with smooth and bounded but non-Lipschitz coefficients, for which they show that the Euler--Maruyama method converges, but also that any numerical method will converge arbitrarily slowly, see \cite{hairer2015,jentzen2016,muellergronbach2016,yaroslavtseva2017,gerencser2017,muellergronbach2018}. While these are academic examples, \cite{hefter2019a} shows a similar result for the well-known Heston model from mathematical finance for some parameter configurations. This model is used in practice; the parameters are determined by the market, so we cannot choose them freely; for some cases simulations in this model will not converge in reasonable time, which would however be essential for pricing financial derivatives.
All these examples show impressively that irregularities can indeed matter for numerics.

However, in the above examples the problems are not about discontinuities. In \cite{goettlich2017} the authors performed a numerical study on SDEs with discontinuous drift coefficient. They observe a strange effect:~the convergence order seems to be too low in the case where the drift points away from the discontinuity. A similar effect was observed by \cite{PS20} in presence of jump noise. A possible reason for this can be found in rare-event simulation; too few paths ``find'' the discontinuity. For the following SDE with outward pointing drift
\begin{align*}
 dX_t =  -(0.5-2 \sign(X_t) )dt + dW_t,\qquad
  X_0=0.1
\end{align*}
the left-hand side plot shows a typical path; the right-hand side plot shows a path that hits the discontinuity. This happens so rarely that the author managed to produce the right-hand side path only by manipulating the data.
\begin{center}
\includegraphics[width=0.45\textwidth]{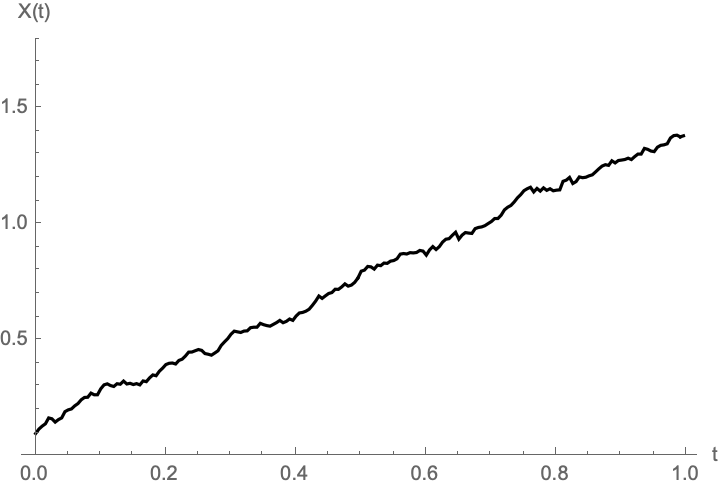}
\includegraphics[width=0.45\textwidth]{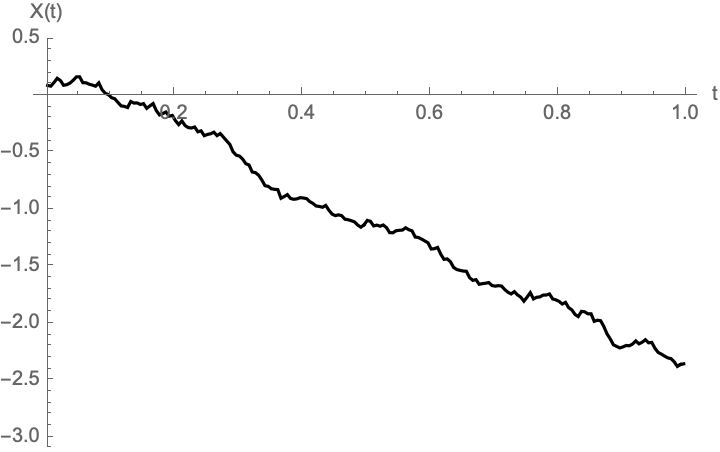}
\end{center}
This shows that in case of a discontinuous drift we have to be even more careful with interpreting numerical results; theoretical results are particularly important.

\section{SDEs with piecewise Lipschitz drift}
\label{pw-lip}

A natural assumption that allows for discontinuities in the drift and that is satisfied in many of the above mentioned applications is piecewise Lipschitz continuity.

\begin{definition}[\text{\cite[Definition 2.1]{LS16}}]\label{def:pwlip-1dim}
A function $f\colon \R\to\R$
is piecewise Lipschitz, if there exist $m\in\N$ and $\zeta_1,\dots,\zeta_m\in
\R$ with $\zeta_1<\dots<\zeta_m$ such that $f$ is Lipschitz on each of the intervals $(-\infty,\zeta_1)$, $(\zeta_1,\zeta_{2})$, \dots, $(\zeta_{m-1},\zeta_{m})$,
$(\zeta_m,\infty)$.
\end{definition}
The points $\zeta_k$ in Definition \ref{def:pwlip-1dim} are possible points of discontinuity of the piecewise Lipschitz function $f$. A multidimensional extension of this definition is available, but requires notions from differential geometry, see \cite{LS17}.

One idea for solving SDEs with piecewise Lipschitz drift coefficient is to use a transformation-based method.
\begin{algorithm}[\text{\citet{LS16,LS17}}]\label{algo-G}
The transformation method works as follows:
\begin{itemize}
\item Construct a transform $G$ depending on the coefficients of the SDE such that the SDE for $G(X)$ has Lipschitz coefficients;
\item compute its inverse $G^{-1}$;
\item define $Z=G(X)$ and calculate the SDE for $Z$ by It\^o's formula;
\item apply the Euler--Maruyama scheme to compute the approximate solution $Z^{\delta}$ to this SDE;
\item define the approximation $\bar X$ of the solution to the original SDE by $\bar X=G^{-1}( Z^{\delta})$.
\end{itemize}
\end{algorithm}
\begin{theorem}[\text{\citet{LS16,LS17}}]\label{conv-trans}
If $\mu$ is piecewise Lipschitz, $\sigma$ is Lipschitz, and $\sigma(\text{points of discontinuity})>0$, then the transformation method has strong convergence order $1/2$.
\end{theorem}
The multidimensional version of Theorem \ref{conv-trans} additionally requires conditions on the hypersurface of discontinuity, see \cite{LS17,LS19}. An overview can be found in \cite{LS17b}.

Now we know a numerical method for SDEs with discontinuous drift of which we know that it converges at the highest rate that can be expected. However, there are two shortcomings of the method. First, in the multidimensional case $G$ has to be inverted numerically, which is costly. Second, we have to know the points of discontinuity for the construction of $G$. Numerically calculating those points, e.g., in the case of applications in optimal control, where the discontinuity is imposed by the control, is also costly. So actually we prefer a simple and explicit method.\\ 

What about the Euler--Maruyama method? We can make use of the transformation idea from above to split the error of the Euler--Maruyama approximation $X^{\delta}$ into the approximation error of the Euler--Maruyama approximation of $Z=G(X)$, that is $\E[   \|Z_t-  Z^{\delta}_t \|^2]$, and the difference between the approximation error of the transformed equation and the transformation of the approximation error of the Euler--Maruyama approximation of $X$, that is $ \E[   \| Z^{\delta}_t-G ( X^{\delta}_t ) \|^2]$:
\begin{align*}
&\E\!\left[   \|X_t- X^{\delta}_t \|^2 \right]
= \E \!\left[ \|G^{-1}(Z_t)-G^{-1}(G ( X^{\delta}_t )) \|^2 \right]\\
&\le 2(L_{G^{-1}}) ^2\E\!\left[   \|Z_t-  Z^{\delta}_t \|^2 \right]
+2(L_{G^{-1}}) ^2 \E \!\left[   \| Z^{\delta}_t-G ( X^{\delta}_t ) \|^2 \right].
\end{align*}
The hard part here is estimating the second error term. We obtain the following result.
\begin{theorem}[\text{\citet{LS18}}]\label{thm-euler}
Under the same assumptions as in Theorem \ref{conv-trans}, the Euler--Maruyama method has essentially strong convergence order $1/4$.
\end{theorem}
So on the one hand we have the transformation method that convergences quickly, but is not simple to implement and potentially costly, and we have the simple and explicit Euler--Maruyama method of which up to now we only have convergence order $1/4$.\\

Can we do better by incorporating more knowledge about the points of discontinuity?
The idea is, instead of using a deterministic grid, to choose the step size $\tau_k$ adaptively in dependence of the current position of the approximation $\appxh$:
\begin{align*}
 \tau_0 = 0, \qquad \tau_{k+1}= \tau_{k} + h(\appxh_{\tau_k},\delta)
 \end{align*}
for a step size function $h$ that decreases whenever $\appxh_{\tau_k}$ gets closer to the discontinuity of the drift and where $\delta$ is again the maximal step size. For the exact choice of $h$, see \cite{NSS19}. The resulting adaptive Euler--Maruyama scheme looks quite similar to the classical one; only the step sizes are now random variables as well:
 \begin{align*}
\appxh_0 =\xi, \qquad  \appxh_{\tau_{k+1}} &= \appxh_{\tau_k}+ \mu(\appxh_{\tau_k})(\tau_{k+1}-\tau_k) + \sigma(\appxh_{\tau_k})(W_{\tau_{k+1}}-W_{\tau_k}).
\end{align*}
When studying adaptive methods, it is not sufficient to estimate the rate of convergence, but we also have to estimate the computational cost. For the classical Euler--Maruyama scheme the computational cost is proportional to the (deterministic) number of steps, which is itself proportional to $1/\delta$. A good rate of convergence does not help, if the number of steps $N(h)$, which depends via $h$ on $\appxh$ and hence is random, explodes.
\begin{theorem}[\text{\citet{NSS19}}]
Under the same assumptions as in Theorem \ref{conv-trans}, the rate of convergence of the adaptive Euler--Maruyama scheme is essentially $1/2$. The expected computational complexity is essentially proportional to $1/\delta$, that is for all $\varepsilon \in (0,1)$ there exists $c\in(0,\infty)$ such that
$$ \E[ N(h) ] \leq c \cdot \delta^{-1+\varepsilon}\,.$$
\end{theorem}
So we have found a numerical method that converges at essentially the highest expectable rate at essentially the same computational cost. This is the best result for the multidimensional case we know so far. In the scalar case however there has been a significant improvement to Theorem \ref{thm-euler}. The setup is the same, but in the analysis of the error term $ \E[   \| Z^{\delta}_t-G ( X^{\delta}_t ) \|^2]$, in particular in a critical point of the analysis, a better estimate has been achieved, yielding the optimal convergence order.
\begin{theorem}[\text{\citet{muellergronbach2020}}]
Under the same assumptions as in Theorem \ref{conv-trans} and $d=1$, the Euler--Maruyama method has strong convergence order $1/2$.
\end{theorem}

There is also a very recent result on higher order methods for SDEs with discontinuous drift. The higher order method that has been studied is a transformation-based Milstein scheme and is obtained by replacing the Euler--Maruyama scheme in Algorithm \ref{algo-G} by the Milstein scheme.
\begin{theorem}[\text{\citet{muellergronbach2019b}}]\label{milstein}
The transformation-based Milstein scheme has strong convergence order $3/4$.\footnote{For the assumptions see \cite{muellergronbach2019b}.}
\end{theorem}
In the regular case the convergence order of the Milstein scheme is $1$. So this is the convergence order we would expect also in our case. In Section \ref{rate34} below we will come back to this issue.

In the case where our SDE is scalar and additionally contains Poisson jumps, that is
$$dX_t=\mu(X_t) dt + \sigma(X_t) dW_t + \rho(X_t) dN_t, \qquad X_0=\xi,$$
the optimal convergence order or the Euler--Maruyama method has been proven:
\begin{theorem}[\text{\citet{PS20}}]
In the case of presence of Poisson jumps and $d=1$ the Euler--Maruyama scheme has strong convergence order $1/2$.\footnote{For the assumptions, in particular for those on the jump coefficient $\rho$, see \cite{PS20}.}
\end{theorem}

\subsection*{Additional literature}
 
Further contributions on the topic are \cite{gyongy1998,halidias2008,ngo2016,ngo2017a,ngo2017b}.
For SDEs with one-sided Lipschitz drift, \cite{gyongy1998} proves an almost sure convergence result.
For discontinuous but monotone drift coefficient and additive noise (i.e.~$\sigma\equiv 1$), \cite{halidias2008} provides strong convergence of the Euler--Maruyama scheme.
In the multidimensional case, for one-sided Lipschitz drift that is an appropriate limit of smooth functions, \cite{ngo2016} proves $L^2$-convergence order $1/4$.
In \cite{ngo2017a} they extend this result to not necessarily one-sided Lipschitz drift functions for scalar SDEs,
and in \cite{ngo2017b} they also allow for certain discontinuous diffusion coefficients.
The last three results require a uniformly non-degenerate diffusion coefficient, which is not required in the results presented above.
For scalar SDEs with additive noise the most general result is \cite{gerencser2020}, where they prove an $L^2$-order of essentially $1/2$ even for drift coefficients which are only  bounded and integrable.

Weak convergence results have been obtained in \cite{kohatsu2017,frikha2018}. In 
\cite{etore2013,etore2014} the authors present an exact simulation algorithm for scalar SDEs with drift that is discontinuous in one point, but differentiable everywhere else.

\section{Dependence of the convergence order on the regularity of the drift}

Now that we know quite a number of results for SDEs with discontinuous drift, let us take a step back and ask a more general question:~how does the convergence order of the Euler--Maruyama method depend on the regularity of the drift coefficient? For this we consider a more simple setup, that is we study scalar SDEs with additive noise:
\begin{align*}
dX_t= \mu(X_t) dt +  dW_t, \qquad X_0=\xi,
\end{align*}
where $\mu\colon\R \rightarrow \R$ is the (not necessarily continuous) drift coefficient.

If $\mu$ is bounded and measurable, \cite{Zvonkin1974} ensures existence and uniqueness of a solution. So we do not have to worry about existence and uniqueness and can directly proceed to numerics.

In \cite{NS20} we provide a novel framework for the error analysis:~we decompose the error into a discretisation error and an error coming from approximating a quadrature problem for Brownian motion.
\begin{theorem}[\text{\citet[Theorem 2.4]{NS20}}]\label{quad}  Assume that $\mu$ is bounded and can be decomposed into a regular part $a \in C_b^2( \mathbb{R},\R)$ and an irregular part $b \in L^1(\mathbb{R},\R)$, that is $\mu=a+b$. Then for all $\varepsilon\in(0,1)$, there exists $c\in(0,\infty)$ such that
$$  \E\!\left[|X_T- X^{\delta}_T|^2\right] \leq  c \cdot  \left( \delta^2 + \mathcal{W}^{1-\varepsilon} \right),$$
where
\begin{align*}
\mathcal{W} &= \E\!\left[ \left| \int_0^T   G'(W_s+\xi)  \left(b(W_s+\xi)-  b(W_\es+\xi)\right) ds \right|^2\right],
\end{align*}
and where $G$ is a Zvonkin-type transform (see \cite{Zvonkin1974,NS20}) for the irregular part of the drift.
\end{theorem}
It remains to analyse the error from approximating the quadrature problem $\mathcal{W}$. For this we choose an appropriate function space for $b$ where the regularity of $b$ is determined by a parameter $\kappa$; the higher $\kappa$ the more regular is $b$. Under the additional assumption that $b$ lies in this function space, we get the following overall approximation result, where the error depends, as desired, on the regularity of the drift.
\begin{theorem}[\text{\citet[Corollary 3.9]{NS20}}]\label{appr-reg}
Under the same assumptions as in Theorem \ref{quad} and under the additional assumption that
there exists  $\kappa \in  (0,1)$ such that 
$$ |b|_{\kappa}:= \left( \int_{\mathbb{R}} \int_{\mathbb{R}} \frac{|b(x)-b(y)|^2}{|x-y|^{2 \kappa  + 1}} dx  dy \right)^{\!1/2}< \infty,$$
i.e.~$b$ belongs to the fractional Sobolev-Slobodeckij space of order $\kappa$,
we have for all $\varepsilon\in(0,1)$ that there exists a  constant $c\in(0,\infty)$ such that the Euler--Maruyama scheme satisfies
$$  \left(\E\!\left[\|X_t -  X^{\delta}_t\|^2\right]\right)^{1/2} \le  c \cdot \delta^{(1+\kappa)/2-\varepsilon},  $$
i.e.~the convergence order is essentially $(1+\kappa)/2$.
\end{theorem}
Note that in the case $a\equiv 0$, by \cite{gerencser2020} the above result also holds for $\kappa = 0$.\\

But which kind of functions belong to a fractional Sobolev-Slobodeckij space of order $\kappa$? This finally brings us back to our example from Section \ref{numSDE}.
\begin{example}[\text{\citet{NS20}}]\label{ex-sign}
Let $\mu(x)=\sign(x)$. A decomposition $\mu(x)=a(x)+b(x)$ which satisfies the assumptions of Theorem \ref{appr-reg} for all $\kappa <1/2$ is
\begin{center}
\includegraphics[scale=0.45]{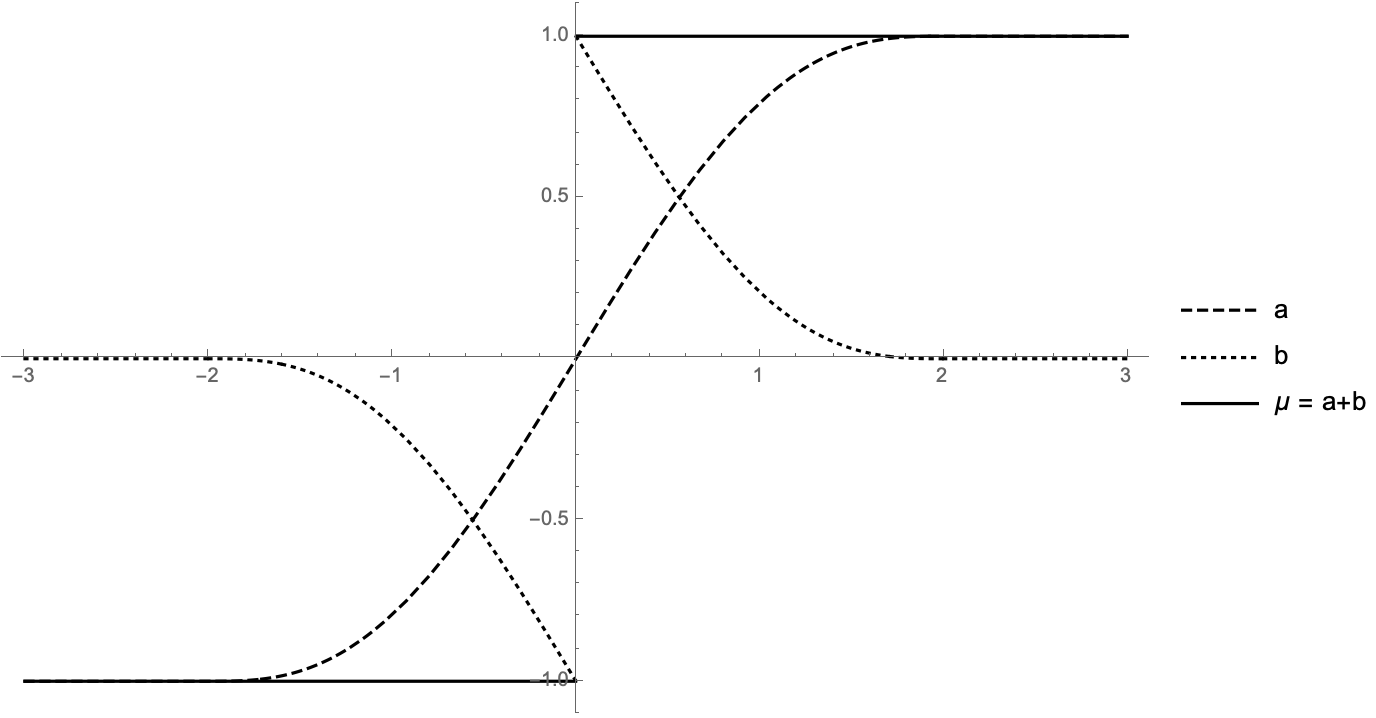}
\end{center}

Hence, the convergence order of the Euler--Maruyama scheme for the SDE with $\mu(x)=\sign(x)$ is essentially $3/4$.
\end{example}

\section{What is behind the rate $3/4$?}
\label{rate34}

Note that for SDEs with additive noise and Lipschitz drift the convergence order of the Euler--Maruyama scheme is 1. In Example \ref{ex-sign} we only obtain rate $3/4$.
Recall that we made a similar observation in Section \ref{pw-lip}:~under classical assumptions the Milstein scheme has convergence order 1, but in the case of a discontinuous drift \cite{muellergronbach2019b} (here Theorem \ref{milstein}) show convergence order $3/4$ of a transformation-based Milstein scheme. Were our estimates too coarse, or is there a structural difference between the convergence orders in the classical case and in the case of a discontinuous drift?
This question is answered by the following very recent result.
\begin{theorem}[\text{\citet{muellergronbach2020b}}]
For scalar SDEs with additive noise the convergence order of any numerical method on a finite deterministic grid is at most $3/4$.\footnote{For the assumptions see \cite{muellergronbach2020b}.}
\end{theorem}
Previous work on lower error bounds for scalar SDEs with discontinuous drift can be found in \cite{hefter2019b}.

The last theorem shows that there is indeed a structural difference between SDEs with Lipschitz drift and SDEs with discontinuous drift. So yes indeed, discontinuities do matter!


%



\end{document}